\documentclass[a4paper,11pt]{amsart}

\usepackage{amsmath,amssymb,amsthm,latexsym,amscd,mathrsfs}
\usepackage{indentfirst}
\usepackage{stmaryrd}
\usepackage{graphicx}
\usepackage{extarrows}
\usepackage{multirow}
\usepackage{latexsym}
\usepackage{amsfonts}
\usepackage{color}
\usepackage{pictexwd,dcpic}
\usepackage{graphicx}
\usepackage{psfrag}
\usepackage{hyperref}
\usepackage{comment}
\usepackage{enumerate}

\numberwithin{equation}{section} \theoremstyle{plain}
\newtheorem{thm}{Theorem}[section]

\newtheorem{prop}[thm]{Proposition}

\makeatletter

\def\<{\langle}
\def\>{\rangle}
\def\({\left(}
\def\){\right)}
\def\[{\left[}
\def\]{\right]}

\makeatother
\title[Complete hypersurfaces in $R^{n+1}$ with constant mean and scalar curvature ]
{Complete hypersurfaces in $R^{n+1}$ with constant mean and scalar curvature }
\author[J.Q. Ge]{Jianquan Ge}
\address{School of Mathematical Sciences, Laboratory of Mathematics and Complex Systems, Beijing Normal
University, Beijing 100875, P.R. CHINA.}
\email{jqge@bnu.edu.cn}

\author[Y. Tao]{Ya Tao$^{*}$}
\address{Chern Institute of Mathematics,
Nankai University, Tianjin 300071, P. R. China.}
\email{tao-ya@mail.nankai.edu.cn}

\date{}
\thanks {$^{*}$ the corresponding author.}
\thanks{J. Q. Ge is partially supported by NSFC (No. 12571049, 12371048),
and the Fundamental Research Funds for the Central Universities.}
\begin{document}
\maketitle

%%%%%%%%%%%%%%%%%%%%%
\begin{abstract}
In this paper, we investigate the rigidity problems of complete hypersurfaces with constant mean curvature and constant scalar curvature in Euclidean spaces. Firstly, under some conditions of Gaussian-Kronecker curvature, we provide  characterizations for the unsolved cases of N\'u\~nez's theorems in dimensions 4 and 5, as well as several rigidity results under some conditions of $r$-th mean curvatures. Moreover, for the case of dimension 6, we also present analogous rigidity results. Finally, for general dimensions, we offer a rigidity theorem under similar pinching conditions.
\end{abstract}
%%%%%%%%%%%%%%%%%%%%%%%

%%%%%%%%%%%%%%%%%%%%%%
\section{Introduction}\label{sec1}
Klotz and Osserman's \cite{TKRO} well-known result states that the only complete surfaces in the Euclidean $3$-space $\mathbb{R}^3$ with non-zero constant mean curvature whose Gaussian curvature does not change sign are the round spheres and the circular cylinders. For higher dimensions, Nomizu and Smyth \cite{KNBS} proved that a complete hypersurface in $\mathbb{R}^{n+1}$ with constant mean curvature, non-negative sectional curvature and constant scalar curvature is a generalized cylinder  $\mathbb{R}^{n-d} \times\mathbb{S}^d$, for some $0\leq d\leq n$.

Cheng and Yau \cite{SYCSTY1, SYCSTY2} proved that the above mentioned result of Nomizu and Smyth holds without the assumption that the scalar curvature is constant, or the mean curvature is constant, respectively. Under the assumption of non-negative sectional curvature, Hartman \cite{PH} established that a complete hypersurface in $\mathbb{R}^{n+1}$ with positive constant $r$-th mean curvature (for some $r=1,...,n$) is $\mathbb{R}^{n-d} \times\mathbb{S}^d$, for some $r\leq d\leq n$.

Motivated by the above discussion, mathematicians naturally seek to generalize the classification theorems by removing the non-negative sectional curvature condition. For $3$-dimensional hypersurfaces, Huang \cite{XGH} partially addressed this problem,  and in 1994, Cheng-Wan \cite{CW} gave a complete classification of constant mean curvature hypersurfaces in $\mathbb{R}^4$ with constant scalar curvature. Compared with Huang's result \cite{XGH}, Cheng-Wan's work \cite{CW} advanced the field in two pivotal aspects: one is proving that the assumption $R\geq 0$ is superfluous, and the other is deriving an explicit description of hypersurfaces when $R=0$. We now state Cheng-Wan's theorem for $H\neq 0$ as follows.
\begin{thm}\label{Thm}$($\cite{CW}$)$
A complete hypersurface of $\mathbb{R}^4$ with non-zero constant mean curvature and constant scalar curvature is one of $\mathbb{R}^{3-d} \times\mathbb{S}^d$, for $d=1,2,3$. 
%\begin{itemize}
%\item[$(1)$] $\mathbb{S}^3(c), c>0$;
%\item[$(2)$] $\mathbb{R}^k\times\mathbb{S}^{3-k}(c), k=1,2$.
%\end{itemize}
\end{thm}
Under some additional assumptions, N\'u\~nez \cite{RAN} extended Cheng-Wan's \cite{CW} theorem to dimensions 4 and 5.
\begin{thm}\label{Thm1}$($\cite[Theorem 1.1]{RAN}$)$
Let $M^4$ be a complete and connected Riemannian manifold with constant scalar curvature $R$, and
$f:M^4\rightarrow \mathbb{R}^5$ be an isometric immersion with non-zero constant mean curvature $H$.
If $R\geq \frac{2}{3}H^2$, then $R=H^2$, $R=\frac{8}{9}H^2$ or $R=\frac{2}{3}H^2$. When $R=H^2$, $ f(M^4)=\mathbb{S}^4(\frac{1}{|H|})$ and when $R=\frac{8}{9}H^2$, $ f(M^4)=\mathbb{R}\times\mathbb{S}^3(\frac{3}{4|H|})$.
\end{thm}

\begin{thm}\label{Thm3}$($\cite[Theorem 1.2]{RAN}$)$
	Let $M^5$ be a complete and connected Riemannian manifold with constant scalar curvature $R$, and $f:M^5\rightarrow\mathbb{R}^6$ be an isometric immersion with non-zero constant mean curvature $H$ and non-negative $4$-th mean curvature $H_4$. If $R\geq \frac{5}{8}H^2$, then $R=H^2$, $R=\frac{15}{16}H^2$, $R=\frac{5}{6}H^2$ or $R=\frac{5}{8}H^2.$ When $R=H^2$, $f(M^5)=\mathbb{S}^5(\frac{1}{|H|})$ and when $R=\frac{15}{16}H^2$, $f(M^5)=\mathbb{R}\times\mathbb{S}^4(\frac{4}{5|H|})$.
\end{thm}

 It is noted that neither of Theorem \ref{Thm1} or Theorem \ref{Thm3} has complete classification under their assumptions. In the first part of this paper, we prove that in the exceptional case $R=\frac{2}{3}H^2$ of Theorem \ref{Thm1}, if the Gauss-Kronecker curvature $H_4\geq 0$,  then $f(M^4) =\mathbb{R}^2\times\mathbb{S}^2(\frac{1}{2|H|})$. And in the case $R=\frac{5}{8}H^2 $ of Theorem \ref{Thm3}, if the Gauss-Kronecker curvature $H_5\leq 0$,  then $f(M^5) =\mathbb{R}^3\times\mathbb{S}^2(\frac{2}{5|H|})$. In addition, under certain conditions involving scalar curvature and $r$-th mean curvature $(r=3,4,5)$, we also obtain some other rigidity results. 

\begin{thm}\label{ThmA}
Let $M^4$ be a complete and connected Riemannian manifold with constant scalar curvature $R$, and
$f:M^4\rightarrow \mathbb{R}^5$ be an isometric immersion with positive constant mean curvature $H$. Then 
\begin{itemize}
	\item[$(1)$] When $R=\frac{2}{3}H^2$, if the  Gauss-Kronecker curvature $H_4\geq 0$, then $f(M^4) =\mathbb{R}^2\times\mathbb{S}^2(\frac{1}{2|H|})$. 
	\item[$(2)$] There is no such hypersurface $M^4$ in $\mathbb{R}^5$ with $R>0$ and $H_4\leq \epsilon<0$.
	\item[$(3)$] There is no such hypersurface $M^4$ in $\mathbb{R}^5$  with $R=\frac{2}{3}H^2$ and $H_3\leq \epsilon<0$.
\end{itemize}
 Here $\epsilon$ is an arbitrarily small negative real number.
\end{thm}

\begin{thm}\label{ThmB}
Let $M^5$ be a complete and connected Riemannian manifold with constant scalar curvature $R$, and let $f:M^5\rightarrow\mathbb{R}^6$ be an isometric immersion with positive constant mean curvature $H$ and non-negative $4$-th mean curvature $H_4$. Then 
\begin{itemize}
\item[$(1)$] When $R=\frac{5}{8}H^2$, if the  Gauss-Kronecker curvature $H_5\leq 0$, then $f(M^5) =\mathbb{R}^3\times\mathbb{S}^2(\frac{2}{5|H|})$.
\item[$(2)$] If $R\geq \frac{5}{8}H^2$ and $H_5\geq \epsilon>0$, then $R=H^2$ and $f(M^5)=\mathbb{S}^5(\frac{1}{|H|})$.
\end{itemize} 
 Here $\epsilon$ is an arbitrarily small positive real number.
\end{thm}

Another major result of this article is for $6$ dimensional hypersurfaces. Upon imposing similar conditions, we establish the following theorems analogous to Theorems \ref{Thm1} and \ref{Thm3}.
\begin{thm}\label{ThmC}
	Let $M^6$ be a complete and connected Riemannian manifold with constant scalar curvature $R$, and let $f: M^6\rightarrow \mathbb{R}^7$ be an isometric immersion with non-zero constant mean curvature, non-negative $4$-th mean curvature $H_4$ and zero $5$-th mean curvature $H_5\equiv 0$. If $R\geq\frac{3}{5}H^2$, then $R=\frac{3}{5}H^2$, $R=\frac{4}{5}H^2$ or $R=\frac{9}{10}H^2$.
\end{thm}

%Analogous to the results of Theorem \ref{ThmA} and Theorem \ref{ThmB}, we present a characterization of the critical case for Theorem \ref{ThmC}. In comparison to the original theorem, this only requires the $5$-th mean curvature to be non-positive.
\begin{thm}\label{ThmC1}
	Let $M^6$ be a complete and connected Riemannian manifold with constant scalar curvature $R$,  and let $f: M^6\rightarrow \mathbb{R}^7$ be an isometric immersion with positive constant mean curvature, non-negative $4$-th mean curvature $H_4$ and non-positive $5$-th mean curvature $H_5$. When $R=\frac{3}{5}H^2$, if the Gauss-Kronecker curvature $H_6\geq 0$, then $f(M^6)=\mathbb{R}^4\times\mathbb{S}^2(\frac{1}{3|H|})$.
\end{thm}

The last part of this paper aims to establish an analogue of the following rigidity theorem of \cite{RAN} with opposite pinching condition for general dimension. 
\begin{thm}\label{Thm2}$($\cite[Theorem 1.3]{RAN}$)$
	Let $M^n$ be a complete and connected Riemannian manifold of dimension $n\geq 3$ and constant scalar curvature $R$, and let $f: M^n\rightarrow \mathbb{R}^{n+1}$ be an isometric immersion with non-zero constant mean curvature $H$. If $HH_3\geq 0$ and $0\leq R\leq\frac{nH^2}{2(n-1)}$, then $R=0$ or $R=\frac{nH^2}{2(n-1)}$. In case that $R=0$, $f(M^n)=\mathbb{R}^{n-1}\times\mathbb{S}^1(\frac{1}{n|H|})$ and in case that $R=\frac{nH^2}{2(n-1)}$, $f(M^n)=\mathbb{R}^{n-2}\times\mathbb{S}^2(\frac{2}{n|H|})$.
\end{thm}
\begin{thm}\label{ThmD}
	Let $M^n$ be a complete and connected Riemannian manifold of dimension $n\geq 3$ and constant scalar curvature $R$. and let $f:M^n\rightarrow\mathbb{R}^{n+1}$ be an isometric immersion with non-zero constant mean curvature $H$. If $R\geq\frac{n(n-2)}{(n-1)^2}H^2$, then $R=H^2$ or $R=\frac{n(n-2)}{(n-1)^2}H^2$. In case that $R=H^2$, $f(M^n)=\mathbb{S}^n(\frac{1}{|H|})$ and in case that $R=\frac{n(n-2)}{(n-1)^2}H^2$, $f(M^n)=\mathbb{R}\times\mathbb{S}^{n-1}(\frac{n-1}{n|H|})$.
\end{thm}

 Note that when $n\geq 3$, the inequality $$\frac{n(n-2)}{(n-1)^2}H^2\geq\frac{nH^2}{2(n-1)}$$ holds, with equality if and only if $n=3$ (in the case of $H\neq0$).  Therefore, by combining the results of Theorem \ref{Thm2} and Theorem \ref{ThmD}, we can recover Theorem \ref{Thm} for $3$ dimensional hypersurfaces with additional conditions $R\geq 0$ and $HH_3\geq 0$.
 
 Note that when $n=4$ or $5$, $R=\frac{n(n-2)}{(n-1)^2}H^2=\frac{8}{9}H^2$ or $\frac{15}{16}H^2$, which aligns with a designated value of $R$ in Theorems \ref{Thm1}, \ref{Thm3}, respectively.

The rest of this paper is organized as follows. In Section \ref{sec2}, we introduce the necessary prerequisites and technical tools required for the subsequent proofs. In Section \ref{sec3}, we provide the proofs of Theorems \ref{ThmA} , \ref{ThmB}, \ref{ThmC}, \ref{ThmC1} and \ref{ThmD}.

\section{Preliminaries}\label{sec2}

%Throughout the paper that $(*)$ is satisfied with equality at $x\in M^n$ means that the inequality at that point is not strict, that is, there exists a unit vector $X\in T_xM$ such that $Ric(X)=\alpha(n,k,H,c).$ If otherwise, we call the inequality $(*)$ strict at $X\in T_xM.$
Let $f:M^n\rightarrow \mathbb{R}^{n+1}$ be an isometric immersion, where $M^n$ is orientable. Let $A$ denote the shape operator of $f$ with respect to a global unit normal vector field $\xi$, and let $\lambda_1,\cdots,\lambda_n$ denote the eigenvalues (principal curvatures) of $A$ with local orthonormal eigenvectors (principal directions) $e_1,\cdots,e_n$. Note that if we label the principal curvatures at each point by the condition $\lambda_1\leq\cdots\leq\lambda_n$, then the principal curvature functions $\lambda_i:M\rightarrow\mathbb{R}$, $i=1,\cdots,n$, are continuous and smooth on an open dense subset of $M^n$.

From the Gauss equation, we can write the average scalar curvature of $M^n$ by
\begin{equation}\label{gauss1}
	R=\frac{2}{n(n-1)}\sum_{i<j}K(e_i,e_j)=\frac{2}{n(n-1)}\sum_{i<j}\lambda_i\lambda_j,
\end{equation}
where $K(e_i,e_j)=\lambda_i\lambda_j$ is the sectional curvature on the plane spanned by $e_i,e_j$. Denote the mean curvature by
 $$H=\frac{1}{n}\mathrm{tr}(A)=\frac{\lambda_1+\cdots+\lambda_n}{n}.$$
 Then we have $|A|^2=\mathrm{tr}(A^2)=\sum_i\lambda_i^2$, 
\begin{equation}\label{guass2}
	n^2H^2=(\sum_{i=1}^{n}\lambda_i)^2=\sum_{i=1}^{n}\lambda_i^2+\sum_{i\neq j}\lambda_i\lambda_j=|A|^2+n(n-1)R.
\end{equation}
%and let $\operatorname{II}: TM\times TM\rightarrow NM$ be the second fundamental form.
Denote the trace-free shape operator $\phi=A-HI$. By direct calculation, $\phi$ satisfies the following three equalities:
\begin{equation}\label{eq2}
	\begin{aligned}
		\mathrm{tr}(\phi)&=0,\\
		|\phi|^2&=\mathrm{tr}(\phi^2)=|A|^2-nH^2,\\
		\mathrm{tr}(A^3)&=\mathrm{tr}(\phi^3)+3H|\phi|^2+nH^3.
	\end{aligned}
\end{equation}

Recall that the $r$-th mean curvature $H_r$ $(1\leq r\leq n)$ is given by: 
\begin{equation}\label{rmean}
	\frac{n!}{r!(n-r)!}H_r=S_r:=\sigma_r(\lambda_1,\cdots,\lambda_n)=\sum_{i_1<\cdots<i_r}\lambda_{i_1}\cdots\lambda_{i_r},
\end{equation}
where $\sigma_r$ is the $r$-th elementary symmetric function on $\mathbb{R}^n$. 
 Notice that $H_1$ and $H_n=\lambda_1\lambda_2\cdots\lambda_n$ are assigned to the mean curvature and the Gauss-Kronecker curvature, respectively. Additionally, the following equality is required in the proof:
\begin{equation}\label{eq}
	\sigma_r(x)=x_i\sigma_{r-1}(\hat{x_i})+\sigma_r(\hat{x_i}), \quad \mbox{for}\ i,r=1,...,n, \quad x\in\mathbb{R}^n,
\end{equation}
where $\sigma_0=1$ and $\sigma_r(\hat{x_i})=\sigma(x_1,...,x_{i-1},x_{i+1},...,x_n)$.

Next, for the completeness of the narrative in the article, we now present the foundational tools employed in our proofs. These tools have been systematically introduced in Section $3$ of \cite{RAN}.
\begin{prop}$($\cite[p.668]{LHZ}\label{Prop}$)$
	%Let $f:M^n\rightarrow \mathbb{R}^{n+1}$ be an isometric immersion, $p\in M^n.$ Denote by $\lambda_1,\cdots,
%	\lambda_n$ the principal curvatures of $M^n$ with orthonormal principal directions $e_1,\cdots,e_n$, and $K_{ij}=\lambda_i\lambda_j$. Then
With notations above, we have	
\begin{equation*}%\label{lap1}
		\frac{1}{2}\Delta|A|^2=|\nabla A|^2+n\sum_i\lambda_i\mathrm{Hess}H(e_i,e_i)+\sum_{i<j}(\lambda_i-\lambda_j)^2\lambda_i\lambda_j.
	\end{equation*}
\end{prop}

%Since $K_{ij}=\lambda_i\lambda_j$,  we rewrite the preceding formula as: 
%	\begin{equation*}%\label{lap2}
%	\frac{1}{2}\Delta|A|^2=|\nabla A|^2+n\sum_i\lambda_i\mathrm{Hess}H(e_i,e_i)+nH\mathrm{tr}(A^3)-|A|^4.
%\end{equation*}
 
In this paper, we consider hypersurfaces with constant mean curvature and constant scalar curvature in Euclidean spaces, thus by (\ref{guass2}), $|A|^2$ is constant and thus all principal curvatures are bounded. Then it follows that
\begin{equation}\label{lap3}
	\begin{aligned}
		0=\frac{1}{2}\Delta|A|^2&=|\nabla A|^2+\sum_{i<j}(\lambda_i-\lambda_j)^2\lambda_i\lambda_j\\
		&\geq\sum_{i<j}(\lambda_i-\lambda_j)^2\lambda_i\lambda_j.
	\end{aligned}
\end{equation}

%Next, we state the lemma that will be used to estimate the laplacian of $|A|^2$, and the estimate will be given and utilized in the proof of Theorem \ref{ThmC}.
%\begin{lem}\label{Lem}$($\cite{MO},\cite{HAMdC}$)$
%	Let $\mu_1,\mu_2,\cdots,\mu_n, n\geq 3$, be real numbers such that $\sum_i\mu_i=0$and $\sum_i\mu _i^2=\beta^2$. Then 
%	\begin{equation}\label{eq3}
%		-\frac{n-2}{\sqrt{(n(n-1)}}\beta^3\leq\sum_i\mu_i^3\leq\frac{n-2}{\sqrt{(n(n-1)}}\beta^3,
%	\end{equation}
%	and equality holds in \ref{Lem} if and only if at least $(n-1)$ of $\mu _i$ are equal.
%\end{lem}

Another pivotal tool is Smyth and Xavier's celebrated principal curvature theorem:
\begin{thm}\label{PCT}$($\cite{SX}$)$ 
	Let $M^n$ be a complete immersed orientable hypersurface in $\mathbb{R}^{n+1}$, which is not a hyperplane. Let $\Lambda\subset\mathbb{R}$ be the set of non-zero values assumed by the eigenvalues of $A$ and let $\Lambda^{\pm}=\Lambda\cap\mathbb{R}^{\pm}$.
	\begin{itemize}
		\item[(i)] If $\Lambda^+$ and $\Lambda^-$ are both nonempty, $\inf\Lambda^+=\sup\Lambda^-=0$.
		\item[(ii)] If $\Lambda^+$ or $\Lambda^-$ is empty then the closure $\overline{\Lambda}$ of $\Lambda$ is connected.
	\end{itemize}
\end{thm}

We conclude this section by presenting the classification results of isoparametric hypersurfaces in Euclidean spaces \cite{GQTY}. Recall that a hypersurface of a real space form $F^{n+1}(c)$ of constant sectional curvature $c$ is called isoparametric if its principal curvatures are constant. When $c=0$, the classification was given by Levi-Civita \cite{TLC} for $n=2$ and by Segre \cite{BS} for higher dimensions. 
\begin{thm}\label{isopara}$($\cite{TLC},\cite{BS}$)$
	If $f:M^n\rightarrow\mathbb{R}^{n+1}$ is a complete isoparametric hypersurface, then $f(M^n)$ is a hyperplane $\mathbb{R}^n$, a round sphere $\mathbb{S}^n(r)$, or a generalized cylinder $\mathbb{R}^{n-k}\times\mathbb{S}^k(r)$, for some $k=1,...,n-1$ and some $r>0$.
\end{thm}

\section{Proof of Theorems}\label{sec3}
In this section, we will prove Theorems \ref{ThmA}, \ref{ThmB}, \ref{ThmC}, \ref{ThmC1} and \ref{ThmD}.

\subsection{Proof of Theorem \ref{ThmA}}
\begin{proof}
Let $\lambda_1\leq\lambda_2\leq\lambda_3\leq\lambda_4$ denote the principal curvatures of the hypersurface $M^4$.
First, at each point of $M^4$, we have $\lambda_4\geq H>0$. We now prove that when $R>0$, for  arbitrary point on the manifold $M$, we have  $\lambda_3>0$. In fact, we have $\inf|\lambda_3|>0$. Otherwise, there is a sequence $\{p_k\}$ in $M$ such that $\lambda_3(p_k) \rightarrow 0$. Since the principal curvature functions are continuous and bounded, by passing to a subsequence, we can assume that $\lambda_i(p_k)\rightarrow \overline{\lambda}_i$, $\forall i$. Thus
\begin{equation}\label{est1}
\overline{\lambda}_1\leq\overline{\lambda}_2\leq\overline{\lambda}_3=0<H\leq\overline{\lambda}_4.
\end{equation}
Since $H$ and $R$ are constant, we have
\begin{equation*}
	\begin{aligned}
		&4H=\overline{\lambda}_1+\overline{\lambda}_2+\overline{\lambda}_4,\\
		&6R=\overline{\lambda}_1\overline{\lambda}_2+\overline{\lambda}_1\overline{\lambda}_4+\overline{\lambda}_2\overline{\lambda}_4.
	\end{aligned}
\end{equation*}
With the aid of (\ref{est1}), we obtain
$$4H\overline{\lambda}_2=6R-\overline{\lambda}_1\overline{\lambda}_4+\overline{\lambda}_2^2\geq6R>0,$$
Thus, $\overline{\lambda}_2>0$, this contradicts (\ref{est1}). Therefore, $\inf|\lambda_3|>0$, i.e., $\sup\lambda_3<0$ or $\inf\lambda_3>0$. If we had $\sup\lambda_3<0$, then $\sup\Lambda^-=\sup\lambda_3<0$, contradicting Theorem \ref{PCT}. Hence, $\lambda_3\geq \inf \lambda_3>0$. 

According to the theorem, we divide the proof into the following three cases.

Case (1): $H_4\geq 0$. Then at each point of $M^4$, there are the following two possibilities: $0\leq\lambda_1\leq\lambda_2\leq\lambda_3\leq\lambda_4$ or $\lambda_1\leq\lambda_2\leq0<\lambda_3\leq\lambda_4$. We assert that $H_3\geq 0$ holds in both scenarios. The first case is obvious, and for the second case, we proceed as follows. 

By (\ref{gauss1}, \ref{rmean}, \ref{eq}), for every $i$, we have 
$$6R=S_2=\lambda_i\sigma_1(\hat{\lambda_i})+\sigma_2(\hat{\lambda_i})=\lambda_i(4H-\lambda_i)+\sigma_2(\hat{\lambda_i}).$$
Setting $i=4$ in the above equality, by the assumption $R=\frac{2}{3}H^2$, we obtain
\begin{equation*}
	\begin{aligned}
\sigma_2(\hat{\lambda_4})&=\lambda_4^2-4H\lambda_4+6R\\
&=(\lambda_4-2H)^2+6R-4H^2\\
&\geq0.
	\end{aligned}
\end{equation*}
Hence, $$4H_3=S_3=\lambda_4\sigma_2(\hat{\lambda_4})+\sigma_3(\hat{\lambda_4})\geq\sigma_3(\hat{\lambda_4})=\lambda_1\lambda_2\lambda_3\geq 0.$$
In summary, $HH_3\geq 0$, and based on the result of Theorem \ref{Thm2}, it can be inferred that $f(M^4)=\mathbb{R}^2\times\mathbb{S}^2(\frac{1}{2|H|})$.

Case (2): $H_4\leq\epsilon<0$. Then at any point of $M^4$, we have $\lambda_1<0<\lambda_2\leq\lambda_3\leq\lambda_4$. In this scenario, based on the principal curvature theorem (Theorem \ref{PCT}), it follows that $\sup\Lambda^-=\sup\lambda_1=0$, and then there exist a sequence $\{p_k\}$ in $M$ such that $\lambda_1(p_k)\rightarrow 0$. Passing to a subsequence, we can assume that $\lambda_i(p_k)\rightarrow\overline{\lambda}_i,\forall i.$ Thus 
$$0=\overline{\lambda}_1\leq\overline{\lambda}_2\leq\overline{\lambda}_3\leq\overline{\lambda}_4$$
and so $$\lim_{k\rightarrow+\infty}H_4(p_k)=\overline{\lambda}_1\overline{\lambda}_2\overline{\lambda}_3\overline{\lambda}_4=0>\epsilon.$$
 which contradicts the condition $H_4\leq\epsilon<0$ on every point of the manifold.
 
 Case (3): $H_3\leq\epsilon<0$. We first claim that $\lambda_1+\lambda_2\geq 0$.  From $R=\frac{2}{3}H^2$ and (\ref{guass2}) we get $|A|^2=8H^2$. i.e.,
 \begin{equation*}
 	\lambda_1^2+\lambda_2^2+\lambda_3^2+\lambda_4^2=\frac{(\lambda_1+\lambda_2+\lambda_3+\lambda_4)^2}{2},
 \end{equation*}
 by simple calculation, we obtain 
 \begin{equation*}
 (\lambda_1-\lambda_2)^2+(\lambda_3-\lambda_4)^2=2(\lambda_1+\lambda_2)(\lambda_3+\lambda_4),
 \end{equation*}
 and thus $\lambda_1+\lambda_2\geq 0$, which proves the claim. Now for the principal curvatures at any point of $M^4$, there are two cases: $\lambda_1\leq0\leq\lambda_2\leq\lambda_3\leq\lambda_4$ or $0\leq\lambda_1\leq\lambda_2\leq\lambda_3\leq\lambda_4$. It is obvious that the second case contradicts the condition $H_3\leq\epsilon<0$. Therefore, the principal curvatures at any point on the manifold $M$ all satisfy the first inequality relation. Thus, we have $\sup\Lambda^-=\sup\lambda_1$.  According to Theorem \ref{PCT}, $\sup\lambda_1=0$, and then there is a sequence $\{p_k\}$ in $M$ such that $\lambda_1(p_k)\rightarrow0$. By passing to a subsequence, we can assume that $\lambda_i(p_k)\rightarrow \overline{\lambda}_i$, $\forall i$. Thus, $\overline{\lambda}_1=0\leq \overline{\lambda}_2\leq \overline{\lambda}_3\leq \overline{\lambda}_4$, which contradicts the assumption $H_3\leq\epsilon<0$.
\end{proof}

\subsection{Proof of Theorem \ref{ThmB}}
\begin{proof}
Let $\lambda_1\leq\lambda_2\leq\lambda_3\leq\lambda_4\leq\lambda_5$ denote the principal curvatures of the hypersurface $M^5$. First, clearly, we have $\lambda_5\geq H>0$. By  (\ref{gauss1}), (\ref{rmean}) and (\ref{eq}), one has for any $i=1,...,6$, 
$$10R=S_2=\lambda_i\sigma_1(\hat{\lambda_i})+\sigma_2(\hat{\lambda_i})=\lambda_i(5H-\lambda_i)+\sigma_2(\hat{\lambda_i}).$$
From the above equality and the hypothesis $R\geq\frac{5}{8}H^2$, we obtain for any $i$,
\begin{equation}\label{eqB1}
	\begin{aligned}
		\sigma_2(\hat{\lambda_i})&=\lambda_i^2-5H\lambda_i+10R\\
		&=(\lambda_i-\frac{5H}{2})^2+10(R-\frac{5}{8}H^2)\\
		&\geq 0.
	\end{aligned}
\end{equation}
 Next, we claim that at each point of $M^5$, we have $\lambda_4>0$ when the condition $R\geq\frac{5}{8}H^2$ and $H_4\geq 0$ hold. In fact, we can prove $\inf|\lambda_4|>0$. Otherwise, there is a sequence $\{p_k\}$ in $M$ such that $\lambda_4(p_k) \rightarrow 0$. By passing to a subsequence, we can assume that $\lambda_i(p_k)\rightarrow \overline{\lambda}_i$, $\forall i$. Thus
 \begin{equation*}
 	\overline{\lambda}_1\leq\overline{\lambda}_2\leq\overline{\lambda}_3\leq\overline{\lambda}_4 =0<H\leq\overline{\lambda}_5.
 \end{equation*}
 Since $H_4\geq 0$, we have $\overline{\lambda}_3=0$ and therefore,
 $$0<10R=\overline{\lambda}_1\overline{\lambda}_2+\overline{\lambda}_1\overline{\lambda}_5+\overline{\lambda}_2\overline{\lambda}_5\leq\overline{\lambda}_1\overline{\lambda}_2,$$
Thus,
 \begin{equation*}
	\overline{\lambda}_1\leq\overline{\lambda}_2<\overline{\lambda}_3=\overline{\lambda}_4 =0<H\leq\overline{\lambda}_5,
\end{equation*}
which implies $\sigma_2(\overline{\lambda}_2, \overline{\lambda}_3,\overline{\lambda}_4,\overline{\lambda}_5)<0$, contradicting (\ref{eqB1}). Hence,  $\inf|\lambda_4|>0$, and similar to the proof of Theorem \ref{ThmA},
$\lambda_4\geq\inf\lambda_4>0$. This proves the claim.

According to the theorem, we divide the proof into the following two cases.

Case (1): $H_5\leq 0$. Since $H_5\leq 0$,  then the principal curvatures at any point of $M^5$ satisfy either $\lambda_1\leq\lambda_2\leq\lambda_3\leq0<\lambda_4\leq\lambda_5$ or $\lambda_1\leq0\leq \lambda_2\leq \lambda_3\leq\lambda_4\leq\lambda_5$. We now prove that $H_3\geq 0$ holds for both cases under consideration. When $\lambda_1\leq\lambda_2\leq\lambda_3\leq0<\lambda_4\leq\lambda_5$. According to (\ref{rmean}) and (\ref{eq}),
\begin{equation*}
	\begin{aligned}
		5H_4=S_4&=\lambda_4\sigma_3(\hat{\lambda_4})+\sigma_4(\hat{\lambda_4})\\
		&=\lambda_4\sigma_3(\hat{\lambda_4})+\lambda_1\lambda_2\lambda_3\lambda_5\\
		&\geq 0.
	\end{aligned}	
\end{equation*} 
Notice that $\lambda_1\lambda_2\lambda_3\lambda_5\leq 0$,  $\lambda_4>0$, we get $\sigma_3(\hat{\lambda_4})\geq 0$. From (\ref{eqB1}), we have $\sigma_2(\hat{\lambda_4})\geq 0$. Once again, based on (\ref{rmean}) and (\ref{eq}), we have $$10H_3=S_3=\lambda_4\sigma_2(\hat{\lambda_4})+\sigma_3(\hat{\lambda_4})\geq0.$$

When $\lambda_1\leq0\leq \lambda_2\leq \lambda_3\leq\lambda_4\leq\lambda_5$, the proof for this case proceeds analogously to that of the first case. 

In summary, we have established that $HH_3\geq 0$ at every point of the manifold. By the conclusion of Theorem \ref{Thm2}, it follows that $f(M^5)=\mathbb{R}^3\times\mathbb{S}^2(\frac{2}{5|H|})$.

Case (2): $H_5\geq\epsilon>0$, then the principal curvatures at any point of $M^5$ satisfy either $0<\lambda_1\leq\lambda_2\leq\lambda_3\leq\lambda_4\leq\lambda_5$ or $\lambda_1\leq\lambda_2< 0<\lambda_3\leq\lambda_4\leq\lambda_5$.  Let $S$ denote the set of points on the manifold where the principal curvatures satisfy the second type of condition, then $\sup_{\{x\in M\}}\Lambda^-=\sup_{\{x\in S\}}\lambda_2$.  According to Theorem \ref{PCT}, $\sup_{\{x\in S\}}\lambda_2=0$,  and then there exists a sequence of points $\{p_k\}\subset S$ such that  $\lambda_2(p_k)\rightarrow 0$. Since the principal curvature functions are bounded, by passing to a subsequence if necessary, we can assume $\lambda_i(p_k)\rightarrow\overline{\lambda}_i,\forall i.$ Thus 
$$\overline{\lambda}_1\leq\overline{\lambda}_2=0\leq\overline{\lambda}_3\leq\overline{\lambda}_4\leq\overline{\lambda}_5.$$
and so $$\lim_{k\rightarrow+\infty}H_5(p_k)=\overline{\lambda}_1\overline{\lambda}_2\overline{\lambda}_3\overline{\lambda}_4\overline{\lambda}_5=0<\epsilon,$$
which contradicts the condition $H_5\geq\epsilon>0$ on any point of $M^5$. Consequently, the principal curvatures at any point of $M^5$ satisfy the first inequality relation.  According to (\ref{lap3}) and (\ref{gauss1}), we have $$\lambda_1=\lambda_2=\lambda_3=\lambda_4=\lambda_5=H,$$ and so 
$$R=H^2.$$ Based on the conclusion of Theorem \ref{Thm3}, we can see that $f(M^5)=\mathbb{S}^5(\frac{1}{|H|})$.

%and from $H_4=\tilde{\lambda}_1\tilde{\lambda}_3\tilde{\lambda}_4\tilde{\lambda}_5\geq 0$, we obtain $\tilde{\lambda}_1=0$. According to (\ref{lap4}), $\tilde{\lambda}_4=\tilde{\lambda}_5=:\mu$, $\lambda\lambda_3=$

\end{proof}

\subsection{Proof of Theorem \ref{ThmC}}
\begin{proof}
	We choose an orientation so that $H>0$. Denoting by $\lambda_1\leq\lambda_2\leq\lambda_3\leq\lambda_4\leq\lambda_5\leq\lambda_6$ the principal curvatures of $M^6$, one has, at each point of $M^6$, 
	\begin{equation*}
		\lambda_6\geq H>0.
	\end{equation*}
By (\ref{gauss1}), (\ref{rmean}) and (\ref{eq}), one has for any $i=1,...,6$, 
$$15R=S_2=\lambda_i\sigma_1(\hat{\lambda_i})+\sigma_2(\hat{\lambda_i})=\lambda_i(6H-\lambda_i)+\sigma_2(\hat{\lambda_i}).$$
From the above equality and the hypothesis $R\geq\frac{3}{5}H^2$, we obtain for any $i$,
\begin{equation}\label{eqC1}
	\begin{aligned}
		\sigma_2(\hat{\lambda_i})&=\lambda_i^2-6H\lambda_i+15R\\
		&=(\lambda_i-3H)^2+15(R-\frac{3}{5}H^2)\\
		&\geq 0.
	\end{aligned}
\end{equation}

Now we claim that $\inf|\lambda_5|>0$. Suppose, by contradiction, $\inf|\lambda_5|=0$, and then there exists a  sequence $\{p_k\}\subset M$ such that $\lambda_5(p_k)\rightarrow 0$. Since the principal curvature functions are bounded, by passing to a subsequence, we can assume that $\lambda_i\rightarrow\overline{\lambda}_i$, thus 
$$
	\overline{\lambda}_1\leq\overline{\lambda}_2\leq\overline{\lambda}_3\leq\overline{\lambda}_4\leq\overline{\lambda}_5=0<H\leq\overline{\lambda}_6.
$$
From $H_5\leq 0 ~(H_5\equiv 0)$, we get $\overline{\lambda}_4=0$, and then from $H_4\geq 0$, $\overline{\lambda}_3=0$. Therefore 
$$0<15R=\overline{\lambda}_1\overline{\lambda}_2+\overline{\lambda}_1\overline{\lambda}_6+\overline{\lambda}_2\overline{\lambda}_6\leq\overline{\lambda}_1\overline{\lambda}_2,$$
Hence
\begin{equation*}
	\overline{\lambda}_1\leq\overline{\lambda}_2<\overline{\lambda}_3=\overline{\lambda}_4=\overline{\lambda}_5=0<H\leq\overline{\lambda}_6.
\end{equation*}
This implies $\sigma_2(\overline{\lambda}_2, \overline{\lambda}_3,\overline{\lambda}_4,\overline{\lambda}_5,\overline{\lambda}_6 )<0$, contradicting (\ref{eqC1}). Thus the claim is proved. Similar to the proof of Theorem \ref{ThmA}, we obtain 
%Since on the one hand 
%\begin{equation}\label{eqC2}
%	lim_{k\rightarrow\infty}S_3(p_k)=\sigma_3(\overline{\lambda}_1,\overline{\lambda}_2,\overline{\lambda}_3,\overline{\lambda}_4,\overline{\lambda}_5,\overline{\lambda}_6)=\overline{\lambda}_1\overline{\lambda}_2\overline{\lambda}_6>0.
%\end{equation}
%On the other hand
%\%begin{equation*}
%	\begin{aligned}
		%lim_{k\rightarrow\infty}S_3(p_k)&=\overline{\lambda}_1\sigma_2(\overline{\lambda}_2,\overline{\lambda}_3,\overline{\lambda}_4,\overline{\lambda}_5,\overline{\lambda}_6)+\sigma_3(\overline{\lambda}_2,\overline{\lambda}_3,\overline{\lambda}_4,\overline{\lambda}_5,\overline{\lambda}_6)\\
	%	&\leq\sigma_3(\overline{\lambda}_2,\overline{\lambda}_3,\overline{\lambda}_4,\overline{\lambda}_5,\overline{\lambda}_6)\\
	%	&=0
%	\end{aligned}
%\end{equation*}
%contradicting (\ref{eqC2}). 
\begin{equation}\label{eqC3}
	\inf \lambda_5>0.
\end{equation}

Next we divide the proof into two cases.

Case (1): there exists a sequence $\{p_k\}$ in $M$ such that $\lambda_4(p_k)\rightarrow 0$. By passing to a subsequence, we can assume that $\lambda_i\rightarrow\overline{\lambda}_i,\forall i$. By (\ref{eqC3}),
$$
	\overline{\lambda}_1\leq\overline{\lambda}_2\leq\overline{\lambda}_3\leq\overline{\lambda}_4=0<\overline{\lambda}_5\leq\overline{\lambda}_6,
$$ and by (\ref{lap3}),
\begin{equation}\label{eqC7}
	0\geq\sum_{i<j}(\overline{\lambda}_i-\overline{\lambda}_j)^2\overline{\lambda}_i\overline{\lambda}_j.
\end{equation}
According to the condition $H_5\equiv 0$, we obtain $\overline{\lambda}_3=0$, and thus
\begin{equation}\label{eqC4}
	\overline{\lambda}_1\leq\overline{\lambda}_2\leq\overline{\lambda}_3=\overline{\lambda}_4=0<\overline{\lambda}_5\leq\overline{\lambda}_6.
\end{equation}
From (\ref{eq}), (\ref{eqC1}) and (\ref{eqC4}), we obtain
\begin{equation}\label{eqC5}
	\begin{aligned}
		\lim_{k\rightarrow\infty}S_3(p_k)&=\overline{\lambda}_6\sigma_2(\overline{\lambda}_1,\overline{\lambda}_2,\overline{\lambda}_3,\overline{\lambda}_4,\overline{\lambda}_5)+\sigma_3(\overline{\lambda}_1,\overline{\lambda}_2,\overline{\lambda}_3,\overline{\lambda}_4,\overline{\lambda}_5)\\
		&\geq\sigma_3(\overline{\lambda}_1,\overline{\lambda}_2,\overline{\lambda}_3,\overline{\lambda}_4,\overline{\lambda}_5)=\overline{\lambda}_1\overline{\lambda}_2\overline{\lambda}_5
		\geq 0
	\end{aligned}
\end{equation}
On the other hand,
\begin{equation}\label{eqC6}
	\begin{aligned}
		\lim_{k\rightarrow\infty}S_3(p_k)&=\overline{\lambda}_2\sigma_2(\overline{\lambda}_1,\overline{\lambda}_3,\overline{\lambda}_4,\overline{\lambda}_5,\overline{\lambda}_6)+\sigma_3(\overline{\lambda}_1,\overline{\lambda}_3,\overline{\lambda}_4,\overline{\lambda}_5,\overline{\lambda}_6)\\
		&\leq\sigma_3(\overline{\lambda}_1,\overline{\lambda}_3,\overline{\lambda}_4,\overline{\lambda}_5,\overline{\lambda}_6)=\overline{\lambda}_1\overline{\lambda}_5\overline{\lambda}_6
		\leq 0
	\end{aligned}
\end{equation}
It follows from (\ref{eqC4}), (\ref{eqC5}) and (\ref{eqC6}) that $\overline{\lambda}_1=0$.  Hence by (\ref{eqC7}) and (\ref{eqC4}), 
$$\overline{\lambda}_1=\overline{\lambda}_2=\overline{\lambda}_3=\overline{\lambda}_4=0\ \mbox{and}\ \overline{\lambda}_5=\overline{\lambda}_6=3H,$$
and so 
$$R=\frac{3}{5}H^2.$$

Case (2): Suppose now that $\inf|\lambda_4|>0$. By (\ref{eqC3}) and Theorem \ref{PCT},  we have
\begin{equation}\label{eqC8}
	\inf \lambda_4>0
\end{equation}

Case (2-1): Assume that there is a sequence $\{p_k\}$ in $M$ such that $\lambda_3(p_k)\rightarrow 0$. By passing to a subsequence, we can assume that $\lambda_i\rightarrow\overline{\lambda}_i, \forall i$. By (\ref{eqC8}),
\begin{equation}\label{eqC9}
	\overline{\lambda}_1\leq\overline{\lambda}_2\leq\overline{\lambda}_3=0<\overline{\lambda}_4\leq\overline{\lambda}_5\leq\overline{\lambda}_6.
\end{equation}
Since $H_5\equiv0$, one has $\overline{\lambda}_2=0$, and since $H_4\geq 0$, one has $\overline{\lambda}_1=0$. It follows from (\ref{eqC9}) and (\ref{eqC7}) that 
$$\overline{\lambda}_1=\overline{\lambda}_2=\overline{\lambda}_3=0\ \mbox{and}\ \overline{\lambda}_4=\overline{\lambda}_5=\overline{\lambda}_6=2H,$$
and so 
$$R=\frac{4}{5}H^2.$$

Case (2-2): Assume now that $\inf|\lambda_3|>0$, By (\ref{eqC8}) and Theorem \ref{PCT}, one has 
\begin{equation}\label{eqC10}
	\inf\lambda_3>0
\end{equation}

Case (2-2-1): Assume that there is a sequence $\{p_k\}$ in $M$ such that $\lambda_2(p_k)\rightarrow 0$. By passing to a subsequence, we can assume that $\lambda_i\rightarrow\overline{\lambda}_i, \forall i$. By (\ref{eqC10}),
\begin{equation}\label{eqC11}
	\overline{\lambda}_1\leq\overline{\lambda}_2=0<\overline{\lambda}_3\leq\overline{\lambda}_4\leq\overline{\lambda}_5\leq\overline{\lambda}_6.
\end{equation}
Since  $H_5\equiv0$, one has $\overline{\lambda}_1=0$. It follows from (\ref{eqC7}) and (\ref{eqC11}) that 
$$\overline{\lambda}_1=\overline{\lambda}_2=0\ \mbox{and}\ \overline{\lambda}_3=\overline{\lambda}_4=\overline{\lambda}_5=\overline{\lambda}_6=\frac{3}{2}H,$$
and so 
$$R=\frac{9}{10}H^2.$$

Case (2-2-2): Assume now that $\inf|\lambda_2|>0$. By (\ref{eqC10}) and Theorem \ref{PCT}, one has 
\begin{equation}\label{eqC12}
	\inf\lambda_2>0.
\end{equation}
We have two possibilities:

(1) There exists some $p\in M$ such that $\lambda_1(p)=0$. By (\ref{lap3}) and (\ref{eqC12}),
$$\lambda_1(p)=0,\lambda_2(p)=\lambda_3(p)=\lambda_4(p)=\lambda_5(p)=\lambda_6(p)=\frac{6}{5}H,$$
which contradicts the condition $H_5\equiv 0$.

(2) $\lambda_1(x)\neq 0, \forall x\in M$. Then we must have $\lambda_1>0$ everywhere. Otherwise, if we had $\lambda_1(q)<0$ at some point $q\in M$, then by continuity $\lambda_1$ would be negative everywhere, and from \eqref{eqC12} we get $\inf\Lambda^+>0$, contradicting  Theorem \ref{PCT}. Then from (\ref{eqC12}), $H_5>0$, which contradicts the condition $H_5\equiv 0$.

In conclusion, we have proved that $R=\frac{3}{5}H^2$, $R=\frac{4}{5}H^2$ or $R=\frac{9}{10}H^2$. 
\end{proof}

\subsection{Proof of Theorem \ref{ThmC1}}
\begin{proof}
	Let $\lambda_1\leq\lambda_2\leq\lambda_3\leq\lambda_4\leq\lambda_5\leq\lambda_6$ denote the principal curvatures of the hypersurface $M^6$. First, according to the proof of Theorem \ref{ThmC}, we can assume at each point of $M^6$, $\lambda_6\geq H>0$ and $\lambda_5>0$ when the conditions $R\geq \frac{3}{5}H^2$, $H_5\leq 0$ and $H_4\geq 0$ hold.
	
	We claim that $H_3\geq 0$, which implies $HH_3\geq 0$. Consequently, in accordance with Theorem \ref{Thm2}, we have $f(M^6)=\mathbb{R}^{4}\times\mathbb{S}^2(\frac{1}{3|H|})$. Therefore, it suffices to prove  $H_3\geq 0$ in the following. 
	Since $H_6\geq 0$,  then the principal curvatures at any point of $M^6$ satisfy either $\lambda_1\leq\lambda_2\leq0\leq  \lambda_3\leq\lambda_4\leq\lambda_5\leq\lambda_6$ or $\lambda_1\leq\lambda_2\leq \lambda_3\leq\lambda_4\leq 0<\lambda_5\leq \lambda_6$. We will only provide a proof for the points where the principal curvatures satisfy the first inequality relation; the proof for the other case is analogous. 
	According to (\ref{rmean}) and (\ref{eq}),
	\begin{equation*}
		\begin{aligned}
			6H_5=S_5&=\lambda_5\sigma_4(\hat{\lambda_5})+\sigma_5(\hat{\lambda_5})\\
			&=\lambda_5\sigma_4(\hat{\lambda_5})+\lambda_1\lambda_2\lambda_3\lambda_4\lambda_6\\
			&\leq 0.
		\end{aligned}	
	\end{equation*} 
	Notice that $\lambda_1\lambda_2\lambda_3\lambda_4\lambda_6\geq0$,  $\lambda_5>0$, we get 
	\begin{equation}\label{C11}
			\sigma_4(\hat{\lambda_5})\leq0.
	\end{equation}
Similarly, we have
		\begin{equation*}
			15H_4=\lambda_5\sigma_3(\hat{\lambda_5})+\sigma_4(\hat{\lambda_5})\geq 0
	\end{equation*} 
	Then it follows from (\ref{C11}) that $ \sigma_3(\hat{\lambda_5})\geq 0$. 	From (\ref{eqC1}), we have $\sigma_2(\hat{\lambda_5})\geq 0$. Once again, based on (\ref{rmean}) and (\ref{eq}), we have $$20H_3=S_3=\lambda_5\sigma_2(\hat{\lambda_5})+\sigma_3(\hat{\lambda_5})\geq0.$$
	The claim is proved and thus the proof is completed.
\end{proof}

\subsection{Proof of Theorem \ref{ThmD}}
\begin{proof}
Since, by hypothesis, $H$ and $R$ are constant, one has by (\ref{guass2}) that $|A|^2$ is also constant. 
When $R\geq\frac{n(n-2)}{(n-1)^2}H^2$, from (\ref{guass2}) and (\ref{eq2}) we obtain $|A|^2\leq\frac{n^2H^2}{(n-1)}$ and $|\phi|^2\leq\frac{nH^2}{(n-1)}$, and so $$0\leq|\phi|\leq\frac{n|H|}{\sqrt{n(n-1)}},$$
Therefore, 
$$-|\phi|^2-\frac{n(n-2)}{\sqrt{n(n-1)}}|H||\phi|+nH^2\geq 0.$$
According to (3.7) in \cite{RAN} , we have 
\begin{equation}
	0=\frac{1}{2}\Delta|A|^2\geq|\nabla A|^2+|\phi|^2[nH^2-\frac{n(n-2)}{\sqrt{n(n-1)}}|H||\phi|-|\phi|^2]\geq 0.
\end{equation}
Thus, $\nabla A\equiv 0$ and $|\phi|=0$ or $|\phi|=\frac{n|H|}{\sqrt{n(n-1)}}$. From $\nabla A\equiv 0$ one concludes that $f$ is isoparametric [p. 254, \cite{PJR}]. Since $n\geq 3$ and $H\neq 0$, it follows from Theorem \ref{isopara} that 
$$f(M^n)=\mathbb{R}^k\times\mathbb{S}^{n-k}(r),$$
for some $r>0$ and some $0\leq k\leq n-1$. From $|\phi|=0$ or $|\phi|=\frac{n|H|}{\sqrt{n(n-1)}}$, we obtain $R=H^2$ or $R=\frac{n(n-2)}{(n-1)^2}H^2$ by simple calculations. Using (\ref{guass2}), it can be easily verified that $R=H^2$ occurs precisely when $k=0$, and $R=\frac{n(n-2)}{(n-1)^2}H^2$ when $k=1$. Hence $f(M^n)=\mathbb{S}^n(\frac{1}{|H|})$ when $R=H^2$, and $f(M^n)=\mathbb{R}\times\mathbb{S}^{n-1}(\frac{n-1}{n|H|})$ when $R=\frac{n(n-2)}{(n-1)^2}H^2$.
\end{proof}

%%%%%%%%%%%%%%%%%%%%%%%
%\begin{acknow}
%\end{acknow}
%%%%%%%%%%%%%%%%%%%%%%%

%%%%%%%%%%%%%%%%%%%%%%%%%%%%%%%%%%%%%%%%%%%%

\end{document}